# Optimizing Electric Multiple Unit Circulation Plan within Maintenance Constraints for High-Speed Railway System

Jian Li[a], Boliang Lin[*b]

(a.Xindu Planning and Management Bureau of Chengdu City, Sichuan 610599, China)
(b.School of Traffic and Transportation, Beijing Jiaotong University, Beijing 100044, China)

**Abstract**：The Electric Multiple Unit (EMU) circulation plan is the foundation of EMU assignment and maintenance planning,and its primary task is to determine the connections of trains in terms of train timetable and maintenance constraints. We study the problem of optimizing EMU circulation plan, and a 0-1 integer programming model is presented on the basis of train connection network design. The model aims at minimizing the total connection time of trains and maximizing the travel mileage of each EMU circulation, and it meets the maintenance constraints from two aspects of mileage and time cycle. To solve a large-scale problem efficiently, we present a heuristic algorithm based on particle swarm optimization algorithm. At last, we conclude the contributions and directions of future research.

**Keywords**：High-Speed Railway, EMU Circulation Plan, Train Connection Network,0-1 Integer Programming, Particle Swarm Optimization Algorithm

## 1. Introduction

In recent years, the high-speed railway industry has been developed rapidly all over the world, especially in China, which contributes to the amount of trains simultaneously running on the high-speed railway.There were more than 3,100 pairs of trains in total for the national railway transportation during the spring festival of 2016 in China, and the amount of trains running on the high-speed railway was more than 1,900 pairs. In order to meet the increasing transportation tasksof trains on the high-speed railway network and make sure that each train could be assigned a suitable EMU, the efficiency of EMU utilization should be improvedwhile the number of EMUs is not increased.This problem isgraduallybecoming the focus of attention for railway departments, such as the EMU depot, however, it is a hard work to handle. We should take the optimization of EMU circulation planning into account particularly for this problem, since the quality ofEMU circulation plan determines the connectionsof trains,builds an important foundation for EMU assignment and maintenance planning, and affectsthe efficiency of high-speed railway operation largely. If we handle this problem well,it could have three aspects of benefits, 1 reduce the number of EMUs in service to decrease the expensive procurement cost of EMUs,2 improve the efficiency of high-speed railway operation and 3provide much more transportation capability for the increasing transportation.

As transport securityis one of the most important issues for the high-speed railway operation, and EMU is one of the most crucialtransportation resources, we have to make sure that each EMU in service should be in good conditionand its transport security is guaranteed. Therefore, EMU should go back to the depot for maintenance before the accumulated mileage or the accumulated time meets the maintenance cycle since the last maintenance. According to the maintenance content and the maintenance cycle, the maintenance projectsfor "CRH" type EMU are divided into five different grades in China, which are called routine maintenance(grade I), special maintenance(II), and advanced maintenance(grade III to V).The maintenance mileage cycle and maintenance time cycle of the grade I maintenance are 4,000 kilometers and 2 days (that is 2880 minutes) respectively, which are both the shortest among the 5 grades maintenance. In general, only thegrade I maintenance is involved in the process of EMU circulation planning in China.In this paper, the maintenance mentioned in the later parts refers to the grade Imaintenanceparticularly, and it is taken as the maintenance constraint for EMU circulation planning.

According tothe practical situation of EMU management in China, EMU circulation is usually definedas the connectionsof trains undertaken in order by an EMU, and the EMU is to be arranged for maintenance after completing undertaking the trains. Therefore,in this paper we define the EMU circulation as a setwhich consists of a series of trains ordered by their departure time.As maintenance must be arranged at the depot to which the EMU belongs,the depot plays a role of both the departure place and arrival place of EMU circulation.That is to say the departure station of the first train and the arrival station of the last train are the same station which is adjacent to the depot. In the following(Fig.1.), we give a more understandable explanation with the schematic diagram.

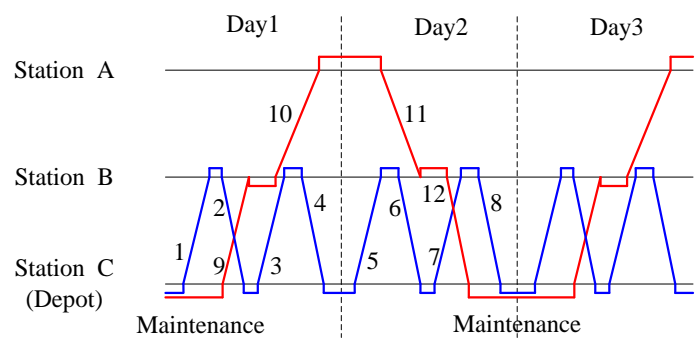

**Fig.1.**Schematic diagram ofEMU circulation

As shown in Fig.1, station C is adjacent to the depot, and if an EMU arrives at station C, it could go back to the EMU depot for maintenance if necessary. EMU could not get maintenance after arriving at Station A and Station B since they are both the turn-back stations for EMU with no depot near them, EMU could not get maintenance after arriving at the two stations, but it could get some service works such as cleaning and parking. In this figure, there are 12 trains in total. The trains with number from 1 to 8 are included in one EMU circulation, and they are sequenced according to the departure time. An EMU undertakes trains 1 to 4one by one on the first dayand comes back to the depot for parking or cleaning in the first night,then itcontinues to undertake the other four trains (5-8) in order on the next day and comes back to the depot for maintenance again in the next night. The trains with number from 9 to 12 consist of another EMU circulation:anEMU departs from station C with undertaking train 9 which is from station C to station B; and then the EMU continues to undertake train 10 from station B to station A;next, the EMU stays at station A for parking orcleaning the first night, on the next day it continues to undertake train 11 and train 12 in sequence, then, it goes back to the depot for the next maintenance.

As the EMU circulation stated above, we concludethat the main task of EMU circulation planning is to determine the connectionsof trainsin terms of the given train timetablewith some objectives pursued and some constraints met in the process.The connectionsbetween any two trains must meet the constraints of connection time and connection station, as well as the limit of the travel mileage and travel timethat each EMU circulation should be within the maximal maintenance cycle. In addition, the connectionsof trains should be as compact as possibleto minimize the total connection time of trainsand to help improve the efficiency of EMU utilization.

EMU circulation plan is usually madeby the management department of EMU with the management department of railway passenger transportationsynergistically, because the EMU assignment and maintenance planning and the crew scheduling are both involved in the problem of EMU circulation planning.However, the crew scheduling problem will not be considered in this paper, and we only study the optimization method for EMU circulation planning from the perspective of the management department of EMU in order to reduce the complexity of this optimization problem. In addition, there are many trains in practicebut not all of them share the same type, and the trains and the EMUs should match each other in the type.What's more, some trains should be assigned coupled EMUs. Besides, the EMUs arrangedto undertake the trains may belong to different depots, that is to say there are more than one EMU depots which are responsible for undertaking the transportation taskstogether. Therefore, it is obvious that EMU circulation planning is a very complex problem with so many situations. However, this complicated problem can be divided into a few simpler sub-problems, in whichthe EMUs included in each sub-problem are of the same type and belong to the same depot. Based on this simplified view, we can tackle each sub-problem by the same optimization method. In this paper, we research the optimization method for EMU circulation planning on the condition that there is only one type for the EMUs and all of thembelong to one depot.Moreover, empty runs are not allowed in this paper, because it occurs infrequently in practice and it is expensive to dispatch for empty runs, so we assume that the trains in the train timetable are in pairs.

The remainder of this paper is organized as follows.In Section 2, we review the literatures related to EMU circulation planning. Section 3 presents the mathematical optimization model of EMU circulation planning as well as the train connection network establishing. Section 4 is about the heuristic algorithm design based on particle swarm optimization algorithm and characteristicof EMU circulation. Conclusion and directions of future research are outlined in Section 5.

## 2. Literature Reviews

Electric Multiple Unit (EMU) or rolling stock is an important transportation resource of railway system, and is of expensive procurement cost. For example, a normal EMU consisting of 8 cars in China is worth about 30,000,000 dollars. The efficient circulation is one of the main objectives for railway departments, which could improve the transportation efficiency and reduce the transportation cost. Therefore, many experts and scholars around the world have focused on this problem or related ones, and a series of research results have been achieved.

As the passenger transportation demand varies a lot during the day, much more transportation capacity is required for the peak time than the off-peak time. In other words, more rolling stocks should be assigned during the peak hours of the day. About this problem, Alfieri et al. (2006) took the seat demand of passengers in peak and off-peak hours into account in their research, and made the train units to be coupled or uncoupled at a certain station, then they presented a solution approach for improving circulation efficiency of rolling stock, which was on the basis of an integer multi-commodity flow model with several additional constraints related to the shunting processes at the stations. Fioole et al. (2006) proposed an extended optimization model on the basis of an existing rolling stock model for routing train units, in their research they not only considered several objective criteria related to operational costs, service quality and reliability of the railway system, but also involved the underway combining and splitting of trains in their extended model. At last, a real-life case study was designed in the two research respectively on the same basis of NS Reizigers and got good results. What's more, the solving algorithm for the problem of rolling stock circulation is as important as the optimization model, and it influences the solutionquality and the solving efficiency largely, a few researchers focus on the solving approach of this problem and get some achievements. Peeters and Kroon (2008) considered the process of adding or removing from the trains at some stations for the train units and described a model for this problem, then they described a branch-and-price algorithm to solve it, the solution generated by the approach was evaluated from the aspects of the service to the passengers, the robustness and the cost of the circulation.

Hong et al. (2009) focused on covering a weekly train timetable with minimal working days of a minimal number of train-sets and presented a two-phased train-set routing algorithm for the problem of the high-speed railway systems in their research, they relaxed maintenance requirements and obtained minimum cost routes at first, and then generated maintenance-feasible routes. At last the algorithm was applied to the Korea high-speed railway system and the solution was an 8.8% improvement compared with a set partitioning approach. The circulation efficiency of rolling stock is influenced by the transportation mode sometimes, especially for the transportation on a railway network rather than a single railway line. Zhao et al. (1997) analyzed the different transportation modes of high-speed railway passenger trains, and then they pointed out that scheduling trains on uncertain railroad region was the better any other modes, and this achievement lays an important foundation for the latter research. Miao et al. (2010) described the problem of EMU circulation as multiple TSP (travelling sales-man problem) with replenishment, and proposed a multi-objective integer programming model based on a given train timetable, then they designed a hierarchical optimization heuristic algorithm based on the connecting network, which sought for minimizing the number of EMUs. Li et al. (2013) analyzed the optimization objectives and primary constraints of EMU circulation plan without considering multi-depot and multi-type of EMUs, based on which they designed the EMU connecting network, and established an integer programming model of EMU circulation plan, then they proposed an improved particle swarm optimization algorithm to solve the model. Wang (2012) researched the optimization of operation and maintenance schemes of EMUs in his doctoral dissertation, for the problem of EMU circulation, the research aimed at minimizing the number of EMUs and reducing the maintenance cost, and the type of EMU, coupled and uncoupled situation and maintenance constraints were all included in the optimization model, then a MAX-M1N Ant System was designed to solve the problem. As it is known that EMUs should go back to the EMU depot for periodical maintenance, and this maintenance process influences the circulation efficiency of EMU, so maintenance constraints of EMU should be taken into account for EMU circulation planning. Zhang et al. (2010) considered the dual constraints of both the kilometrage and time for scheduled inspection and repair in their research of EMU assignment, and the computational result showed that the complexity of EMU scheduling was increased evidently, and the number of EMUs in service and the number of inspections and repairs were both increased compared with the single constraint of the kilometrage. Canca et al. (2014) focused on the problem of rolling stock circulation for railway Rapid Transit Systems, and developed a mixed integer programming model for rolling stock circulation plans with considering a rotating maintenance scheme. The approach could minimize the train empty movements and the train units reserve and balance the workload of the maintenance operation. Giacco et al. (2014) researched the problem of rolling stock rostering and maintenance scheduling for a short-term planning, they viewed the rostering solution as a minimal cost Hamiltonian cycle in a graph with service pairings, empty runs, and short-term maintenance tasks, and a mixed-integer linear-programming formulation for this problem was proposed in their research, then they adopted a commercial MILP solver to compute efficient solutions in a short time. Thorlacius et al. (2015) focused on the large number of practical and railway specific requirements of the Copenhagen suburban passenger railway, and proposed an integrated model for rolling stock planning, which could simultaneously consider all practical requirements for rolling stock planning at DSB S-tog, the suburban passenger train operator of the City of Copenhagen, and a hill climbing heuristic was used for solving the model. Sun at el. (2014) researched the train routing problem of combined with train scheduling on a high-speed railway network from another perspective, and put forward a multi-objective optimization model for train routing, and an improved GA was designed to solve the problem of train routing. The approach took the energy consumption and the user satisfaction into account, this approachwas suitable for a small-scale network. Lai et al. (2015) analyzed the manual process of rolling stock assignment and maintenance planning, and pointed out that it was time consuming and it was impossible to guarantee an optimal solution. Therefore, they developed an exact optimization model and a hybrid heuristic process to improve solution quality and efficiency, and the empirical results indicated that the intelligent approach was much better than the manual process. The authors (2016)of this paper have researched the related problem about motor train set assignment and maintenance scheduling, then a 0-1 integer programming model and a heuristic solution strategy based on particle swarm optimization were proposed, which laid an important foundation for the research in this paper.

The circulation problem of EMU or rolling stock is a hot issue for the railway system, and it is researched widely and deeply from different perspectives with diverse objectives and constraints. In addition, there are a lot of other research fields which are analogous to the EMU or rolling stock circulation, such as the locomotive routing problem, the vehicle routing problem and the aircraft routing problem etc.Cordeau et al. (2000, 2001) researched the problem of locomotive and car assignment for the railway passenger transportation system. According to the analysis of optimization objectives and constraints, they proposed a multi-commodity network flow-based model for assigning locomotives and cars to trains, a branch-and-bound method and a Benders decomposition approach were applied to solve the assigning problem on different solving levels respectively. Lingayaet al. (2002) focused on the car assignment problem of VIA Rail Canada and described an optimization model with meeting the maintenance requirements and minimum connection times. Then they described a heuristic algorithm based on a branch-and-bound method and a column generation method, and a software system based on this approach was developed. Rouillon et al. (2006) focused on the locomotive assignment problem of the railway freight transportation system, and presented an efficient backtracking mechanism that can be added to the heuristic branch-and-price approach while meeting the locomotive availability and maintenance requirements and minimizing the operation cost. Noori et al. (2012) analyzed the different degrees of priority for servicing of trains and took it into consideration while determining locomotive assignment, and a two-phase approach based on genetic algorithm was designed to optimization problem. Ceder (2011) focused on the timetable development and vehicle-scheduling with different vehicles types of the he public-transport operation planning process, and this problem was formulated as a cost-flow network problem with an NP-hard complexity level in the research, then a heuristic algorithm was designed to solve this vehicle scheduling problem. Nikolić and Teodorović (2015) considered the unexpectedly high demand in some nodes of the regular transportation network, and proposed the Bee Colony Optimization algorithm and principles of lexicographic optimization to solve the problem on the basis of a mathematical formulation. During the solving process of the vehicle routing problem, they aimed at minimize the negative consequences of the disturbances caused by the increased customer demand as soon as possible. Samà at el. (2014) researched the real-time problem of aircraft scheduling and routing in terminal control area, and they adopted a simulation method to tackle the problem of aircraft rescheduling and rerouting when the take-off and landing of some aircrafts were disturbed, which aimed at minimizing the delay propagation and reducing the aircraft travel time and energy consumption. Salazar-González (2014) researched a solution approach for an integrated fleet-assignment, aircraft-routing and crew-pairing problem covering the flights of a single day, and this combinatorial problem was considered as a 2-depot vehicle routing problem with driver changes. Then an integer programming model and a heuristic algorithm were described to solve the problem in this research. Dauzère-Pérès at el. (2015) proposed a Lagrangianheuristic framework to solve an integrated problem including rolling stock units and train drivers at the same time, which could handle the weakness of the sequential approach for railway transportation resources planning effectively. All of these research in other fields have something in common with the circulation problem of EMU or rolling stock and this could provide a reference for the research in this paper.

Except for the related literatures mentioned above, there are a lot of other studies which are helpful to our research. Kara and Bektas (2006) extended the classical multiple traveling salesman problem and proposed integer linear programming formulations for single and multi-depot cases. Wang (2014) researched the modeling methods for the traveling salesman problem based on the software Lingo, and proposed solutions for avoiding the separation phenomenon of TSP from three different aspects. These could help us establish a train connection network without any sub-close loop circuits, and the advantages can be absorbed into the model of EMU circulation planning. Ho at el. (2012) applied the particle swarm optimization algorithm to solve the train service timetabling problem. Zhang at el. (2015) surveyed and concluded the application of particle swarm optimization algorithm in different research fields. These research about particle swarm optimization algorithm and the book about this algorithm written by Liu and Niu (2013) could provide references for us to design the solving algorithm in this paper.

In conclusion, most of the relativeresearch for EMU circulation planning at present focus on minimizing the total connection time of trains and reducing the number of EMUs in service for a given train timetable. The researcherstook a series of constraints into account, such as the maintenance cycle constraints of mileage and time, the coupled and uncoupledsituations ofEMU, the matching relation between the EMU and the train, and the empty movement of EMU, and so on. On the basis of existing research achievements, we refer to the idea of establishing a train connection network, and study the optimization method for EMU circulation planning with aiming at minimizing the total connection time of trains and maximizing the travel mileage of each EMU circulation. In this process, themaintenance constraints are taken into account particularly.

## 3. Mathematical Optimization Model

### 3.1 Train connection network

At first, we describe the basic attributes of trains with mathematical notations. It is defined that $V = \{v_i \mid i = 1, 2, \cdots, n\}$ as the set of trains in the train timetable, $v_i$ as the $i$ th train in set $V$, $n$ as the amount of trains in set $V$. For anytrain $v_i$, it has a few attributes including departure station $s_i^d$, departure time $t_i^d$, arrival station $s_i^a$ and arrival time $t_i^a$. Except for the four elements, train $v_i$ also has a travel mileage and a travel time, so we let $d_i^0$ represent the travel mileage (unit: km), let $t_i^0$ represent the traveltime (unit: min). We use $S_{\text{station}}$ torepresent the set of all stations included in the train timetable, so there are inclusion relations $s_i^d \in S_{\text{station}}$ and $s_i^a \in S_{\text{station}}$. We let $S_{\text{maint}}$ represent the set of stations which are adjacent to the depot, and EMU could go to the nearby depot for maintenance if necessary after arriving at these stations, so there is inclusion relation $S_{\text{maint}} \subseteq S_{\text{station}}$. As we research the problem of EMU circulation planning with single-type and single-depot of EMU in this paper, there is only one station in set $S_{\text{maint}}$.

In order to establish a train connection network which reflects EMU circulations, we take the trains in set $V$ as the nodes of the train connection network, and take the connection between each two trains as the arc of the network, which is called train connection arc.Moreover,we add the maintenance into the train connection network in the form of network arc, and it is called EMU maintenance arc.According to the settings of train nodes, train connection arcs and EMU maintenance arcs, it is clear that each EMU circulation could be described as a closed loop circuit, which consists of trains connecting in order and EMU maintenance.In light of the trains and EMU circulations shown in Fig.1, we can draw a train connection network with two closed loop circuits in Fig.2(a).In order to reduce the model complexity of EMU circulation planning, we integrate the two closed loop circuits into one closed loop circuit by connecting them through the EMU maintenance arcs, since the departure place and the arrival place of EMU circulation are both the depot. The single closed loop circuit generated at lastis shown in Fig.2(b), and it includes all the train nodes, train connection arcs, EMU maintenance arcs. Based on the idea of establishinga single closed loop circuit, the first work of EMU circulation planning is to determine the connections of trains and to arrange the maintenancereasonably, so a train connection network withonly oneclosed loop circuitcan be generated. And then, we break down the single closed loop circuit from the EMU maintenance arcs, and a few segments with trains connecting in order are generated, these segments are the EMU circulations.

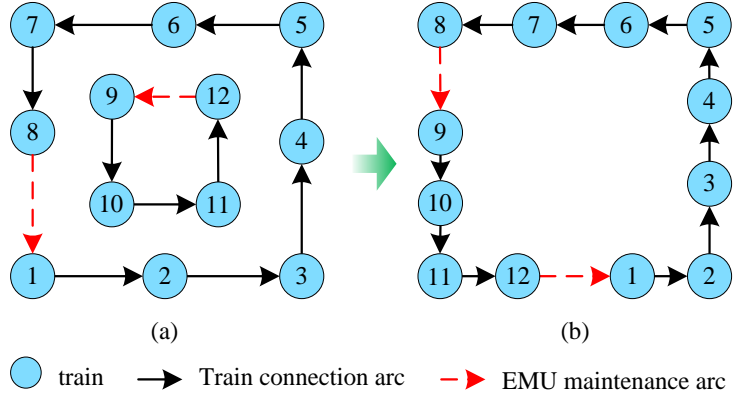

(a)　　　　　　　　　　　(b)

○ train　　→ Train connection arc　　--▶ EMU maintenance arc

**Fig.2.** Schematic diagram of train connection network

On the basis of train connection network establishment, we define a decision variable $x_{ij}$, which represents the connection between train $v_i$ and train $v_j$, if train $v_i$ and train $v_j$ connect together, let $x_{ij}=1$, otherwise let $x_{ij}=0$. Another decision variable $y_{ij}$ is defined to represent whether a maintenance is arranged between train $v_i$ and train $v_j$ if these two trains connect together, if the answer is true, let $y_{ij}=1$, otherwise let $y_{ij}=0$. The value of the two decision variables $x_{ij}$ and $y_{ij}$ are described in formula (1) and formula (2) respectively.

$$x_{ij}=\begin{cases}1,&\text{train }v_i\text{ connects to }v_j\\0,&\text{otherwise}\end{cases} \tag{1}$$

$$y_{ij}=\begin{cases}1,&\text{an maintenance is arranged between train }v_i\text{ and train }v_j\\0,&\text{otherwise}\end{cases} \tag{2}$$

The connection of any two trains should meet a few constraints. As empty runs are not allowed in this paper, the connection between train $v_i$ and train $v_j$ is feasible only if $v_j$ departs form the arrival station of $v_i$. Therefore, we define a parameter $t_{ij}$ which represents the connection time between train $v_i$ and train $v_j$, and it is taken as the generalized cost for the connection of train $v_i$ and train $v_j$. If the arrival station of train $v_i$ and the departure station of train $v_j$ are not the same one, we set the value of connection time $t_{ij}$ as infinite. The value of connection time $t_{ij}$ is determined by formula (3).

$$t_{ij}=\begin{cases}t_j^d-t_i^a&s_i^a=s_j^d,t_j^d-t_i^a\geq t_{connect}\\t_j^d-t_i^a+1440&s_i^a=s_j^d,t_j^d-t_i^a<t_{connect}\\+\infty&\text{otherwise}\end{cases} \tag{3}$$

In formula (3), the parameter $t_{connect}$ represents the minimum time standard of train connection, which mainly includes the passenger service time and the EMU cleaning time. If the calculated value of train connection time $t_j^d-t_i^a$ is less than $t_{connect}$, it indicates that train $v_i$ could not connect to train $v_j$ directly, so we set it that train $v_i$ could connect to a train which is of the same attributes with train $v_j$ in the following day, and the connection time $t_{ij}$ is added by 1440 minutes on the basis of $t_j^d - t_i^a$.

If there is a maintenance arranged between train $v_i$ and train $v_j$, it indicates that train $v_i$ is the last train of the former EMU circulation, and train $v_j$ is the first train of the latter EMU circulation. That is to say another EMU circulation is to be determined after the maintenance, and since EMU circulations have no direct relationship between each other, there is no need to consider the maintenance time and the connection time between train $v_i$ and train $v_j$. So if the arrival station of train $v_i$ and the departure station of train $v_j$ are the same station, and they are adjacent to the depot, which directs to $s_i^d = s_j^a \in S_{\text{maint}}$, we can arrange maintenance according to the maintenance cycle constraints. Because of this, we define an auxiliary parameter $\theta_{ij}$ for the maintenance arrangement, the value of the auxiliary decision parameter $\theta_{ij}$ is determined by formula (4).

$$\theta_{ij} = \begin{cases} 1 & s_i^a = s_j^d \in S_{\text{maint}} \\ 0 & \text{otherwise} \end{cases} \quad (4)$$

For the auxiliary parameter $\theta_{ij}$, if the arrival station of train $v_i$ and the departure station of train $v_j$ are the same station, and they are adjacent to the depot, we let $\theta_{ij} = 1$, otherwise we let $\theta_{ij} = 0$. The matrix of connection time $t_{ij}$ and the matrix of auxiliary parameter $\theta_{ij}$ can both be generated according to the basic information of trains.

### 3.2 Optimization modeling

Based on the analysis of EMU circulation planning problem and the train connection network establishment, we try to present a mathematical optimization model for EMU circulation planning. Except for the variables, parameters and sets defined above, the others involved in the model are defined as follows:

| | |
|---|---|
| $l_j$ : | It represents the accumulated mileage after connecting to train $v_j$ in the train connection network (unit: km). |
| $t_j$ : | It represents the accumulated time after connecting to train $v_j$ in the train connection network (unit: min). |
| $L_{\text{cycle}}$ : | It represents the maintenance mileage cycle (unit: km). |
| $T_{\text{cycle}}$ : | It represents the maintenance time cycle (unit: min). |
| $\lambda$ : | It represents the percentage that the accumulated mileage or time are allowed to exceed the maintenance cycle limit. |
| $\omega_1$ : | It represents the weight coefficient for the total connection time of trains in the optimization objective function. |
| $\omega_2$ : | It represents the weight coefficient for the total difference between the maximal maintenance mileage limit and the accumulated mileage when maintenance occurs. |

In the train connection network, the accumulated mileage $l_j$ and the accumulated time $t_j$ are both variables, and the value of them can be calculated according to the train connections. If train $v_i$ connects to train $v_j$ and there is no maintenance between train $v_i$ and train $v_j$, the value of accumulated mileage $l_j$ is equal to the sum of the former accumulated mileage $l_i$ and the travel mileage $d_j^0$ of train $v_j$. If there is a maintenance between train $v_i$ and train $v_j$, the value of accumulated mileage $l_j$ only equals to the travel mileage $d_j^0$ of train $v_j$ without considering the accumulated mileage $l_i$ before the maintenance. Then the calculation formula of the accumulated mileage $l_j$ is described as follow:

$$l_j = d_j^0 + \sum_{i=1, i \neq j}^{n} x_{ij}(1 - y_{ij})l_i \qquad (5)$$

In the same way, we can get the value of the accumulated time $t_j$. If train $v_i$ connects to train $v_j$ and there is no maintenance between train $v_i$ and train $v_j$, the value of accumulated time $t_j$ is equal to the sum of the former accumulated time $t_i$, the travel time $t_j^0$ of train $v_j$ and the connection time $t_{ij}$. If there is a maintenance between train $v_i$ and train $v_j$, the value of accumulated time $t_j$ equals to the travel time $t_j^0$ of train $v_j$ without considering the accumulated time $t_i$ before the maintenance and the connection time $t_{ij}$. Then the calculation formula of the accumulated time $t_j$ is described as follow:

$$t_j = t_j^0 + \sum_{i=1, i \neq j}^{n} x_{ij}(1 - y_{ij})(t_i + t_{ij}) \qquad (6)$$

Based on the analysis of EMU circulation planning problem and the definition of variables, parameters and sets, we present a 0-1 integer programming model for optimizing EMU circulation plan, and the optimization model is described as follow:

$$\min = \omega_1 \sum_{i=1}^{n} \sum_{j=1, i \neq j}^{n} t_{ij} x_{ij}(1 - y_{ij}) + \omega_2 \sum_{i=1}^{n} \sum_{j=1, i \neq j}^{n} y_{ij}((1 + \lambda)L_{cycle} - l_i) \qquad (7)$$

Subject to:

$$\sum_{i=1, i \neq j}^{n} x_{ij} = 1 \quad j = 1, 2, \cdots, n \qquad (8)$$

$$\sum_{j=1, i \neq j}^{n} x_{ij} = 1 \quad i = 1, 2, \cdots, n \qquad (9)$$

$$y_{ij} \leq x_{ij} \theta_{ij} \quad i, j = 1, 2, \cdots, n \quad i \neq j \qquad (10)$$

$$l_j \leq (1 + \lambda)L_{cycle} \quad j = 1, 2, \cdots, n \qquad (11)$$

$$t_j \leq (1 + \lambda)T_{cycle} \quad j = 1, 2, \cdots, n \qquad (12)$$

$$u_i - u_j + nx_{ij} \leq n - 1 \quad i = 1, 2, \cdots, n \quad j = 2, \cdots, n \qquad (13)$$

$$0 \leq u_i \leq n - 1 \quad i = 1, 2, \cdots, n \qquad (14)$$

$$x_{ij}, y_{ij} \in \{0, 1\} \quad i, j = 1, 2, \cdots, n \quad i \neq j \qquad (15)$$

For this mathematical optimization model, formula (7) is the optimization objective function, and formulas (8) ~ (15) are the constraints. The optimization objective function consists of two parts, the first part is about minimizing the total connection time between any two trains connected by the train connection arcs in the train connection network, and it does not include the connection time between any two trains connected by the EMU maintenance arcs. Seeking for the optimal train connections with the minimum connection time of trains, on one hand, it could help to reduce the amount of EMUs in service. On the other hand, it could help to force the connection of trains included in each EMU circulation much more tightly, which provides more spare time for maintenance arrangement. The second part of the optimization objective function is about minimizing the total difference between the maximal maintenance mileage limit and the accumulated mileage when maintenance occurs, which is the same as maximizing the travel mileage of each EMU circulation without exceeding the maximal maintenance mileage limit. The second optimization objective is conducive to improve the EMU utilization efficiency as much as possible, and it is helpful to reduce the amount of EMUs in service at the same time. In order to reduce the complexity of the model, we define two weight coefficients $\omega_1$ and $\omega_2$ for the two optimization objectives, and transform the multi-objective optimization problem into a single objective optimization problem. For the application of this model, the value of these two weight coefficients $\omega_1$ and $\omega_2$ can be determined according to the practical situation.

For a certain train in the train connection network, there is one and only one train in front which connects to it, and there is one and only one train in behind to which it connects, then these constraints are described in formula (8) and formula (9) respectively. If train $v_i$ connects to train $v_j$, it is possible to arrange maintenance between the two trains only if the maintenance constraints are met, this constraint for maintenance arrangement is described in formula (10). Formula (11) and formula (12) are both the EMU maintenance cycle constraints, in the train connection network they represent that the accumulated mileage and the accumulated time after connecting a train should not exceed the maximal maintenance mileage limit and the maximal maintenance time limit respectively. However, for the practical transportation and management of EMU, the accumulated mileage and the accumulated time after the last maintenance are allowed to exceed the corresponding maintenance cycle within a certain limit, and the percentage for exceeding is less than 10% in general, this practical problem is involved in the description of formula (11) and formula (12). Moreover, in order to make sure that the train connection network generated at last is a single closed loop circuit and to avoid other sub-closed loop circuits being generated, we refer to a classical constraint in the Travelling Salesman Problem which avoids sub-closed loop circuits, then this constraint is described in formula (13) and formula (14). The variables $u_i$ and $u_j$ represent the position number of train $v_i$ and train $v_j$ respectively in the train connection network, and each train has one and only one position number. The last constraint formula (15) represents that the decision variables $x_{ij}$ and $y_{ij}$ should be a 0-1 integer.

## 4. Solution Algorithm

The particle swarm optimization algorithm(PSO) is one of the heuristic algorithms, and it is applied to various research fields since it was proposed in 1990s. The relative researches (Zhang at el., 2015; Niu, 2013) about the particle swarm optimization algorithm indicate that it has a few advantages, such as rapid searching rate, high precision and easy operation, etc. Therefore, in this paper we take the characteristic of EMU circulation into account and design a heuristic algorithm for solving the model on the basis of particle swarm optimization algorithm.

### 4.1 Application strategy of PSO

For the particle swarm optimization algorithm, each particle represents a solution of the optimization problem, so what we should do is to design the form of solution which reflects in each particle. There are two main decision variables for the mathematical optimization model, one is the decision variable of train connection $x_{ij}$ and the other is the decision variable of maintenance arrangement $y_{ij}$. However, it is difficult to design the form of solution by the position of a particle, which could reflect the decision variables $x_{ij}$ and $y_{ij}$. Therefore, in order to tackle this problem, we take advantage of the characteristics of train connection network shown in Fig.2 (b), and we choose a certain node as the first position of the train connection network, then the positions of other nodes are increased progressively according the connecting order. So the amount of positions in the train connection network is the same as the amount of trains, that is to say each train corresponds to a position. On this basis, we define a variable described as $x_d$ for each position in the network, and the parameter $d$ represents the index of positions, $d = 1, 2, \cdots, n$. The variable $x_d$ represents the train assigned to the $d$ th position, and the value of $x_d$ is the train number. Moreover, we define two variables $l_d$ and $t_d$ which represent the accumulated mileage and the accumulated time respectively after assigning a train to the $d$ th position in the train connection network, and the value of $l_d$ and $t_d$ can be determined by referring to the formula (5) and formula (6) respectively. Thus, we set that each dimension of the particle corresponds to each position of the train connection network, and the amount of dimensions of the particle is equal to the amount of positions of the train connection network, that is $n$, and we let $d$ represent the index of dimension of the particle. In addition, we set it that the value of position on each dimension of the particle is equal to the train number which is assigned to the corresponding position of the train connection network. The scale of particle swarm is $N_m$ and the parameter $m$ represents the index of particle, the position vector of each particle is described as $\mathbf{x}_m = (x_{m1}, x_{m2}, \cdots, x_{mn})$, each particle has a velocity vector with the same dimension and it is described as $\mathbf{v}_m = (v_{m1}, v_{m2}, \cdots, v_{mn})$. During the iterative process of optimization computation, there is a personal historical best solution for each particle and there is also a global historical best solution for the particle swarm, which are described as $\mathbf{p}_m = (p_{m1}, p_{m2}, \cdots, p_{mn})$ and $\mathbf{p}_g = (p_{g1}, p_{g2}, \cdots, p_{gn})$ respectively. During the iterative process, the velocity in each dimension of the particle is updated according to the formula (17) descripted as follow:

$$x_{md}^{k+1} = \left\langle x_{md}^k + v_{md}^{k+1} \right\rangle \tag{17}$$

In formula (17), the calculated symbols $<\ >$ means rounding the value to the nearest integer in this paper. The $k$ is index of iterations, and if we let $k$ represent the current iteration, $k+1$ represents the next iteration. $x_{md}^k$ represents the position value in the $d$ th dimension of the $m$ th particle and the value of $x_{md}^k$ is in the value interval $[x_{\min}, x_{\max}]$. If the calculated value gotten by formula (17) is outside of the value interval, we set the value of $x_{md}^k$ as the corresponding boundary value of the value interval. The value interval is set as $[1, n]$, since the value of $x_{md}^k$ is the train number which is an integer. The variable $v_{md}^k$ represents the corresponding velocity of the particle, and in order to avoid falling into the local optimum, we adopt the particle swarm optimization algorithm with inertia weight, and the velocity of the particle is updated by formula (18) which is described as follow:

$$v_{md}^{k+1} = w v_{md}^k + c_1 r_1 (p_{gd} - x_{md}^k) + c_2 r_2 (p_{md} - x_{md}^k) \tag{18}$$

The velocity of the particle also has a value interval $[v_{\min}, v_{\max}]$, and if the calculated value is outside of the value interval, the value of $v_{md}^k$ is set as the corresponding boundary value of the value interval. The parameters $r_1$ and $r_2$ are used for keeping the diversity of the particle swarm, the value of them are both a random number in the value interval $[0,1]$, and they are generated randomly during each iterative process. The parameters $c_1$ and $c_2$ are called learning factor, they make the particle absorb the advantages of the personal historical best solution of itself and the global historical best solution of the particle swarm. In other words, it means that the solution of each iteration approaches to the personal historical best solution and the global historical best solution step by step. The parameter $w$ represents the inertia weight coefficient of velocity, and we adopt a linear function for updating the value of $w$, which decreases while the iteration times increases. Then the calculation formula is described as follow:

$$w = w_{max} - \frac{w_{max} - w_{min}}{k_{max}} \times k \tag{19}$$

In formula (19), the parameters $w_{max}$ and $w_{min}$ represent the maximum value of the inertia weight and the minimum value of the inertia weight respectively, and the parameter $k_{max}$ represents the maximum number of iterations.

As the particle swarm optimization algorithm has a certain randomness, if we solve the model with meeting the constraints strictly, it is possible to get some solutions which are unavailable. Therefore, in order to seek for optimization solution in greater scope and to decrease the complexity, we relax the constraints of maintenance cycle limit, and add a penalty function for the problem of exceeding the maximal maintenance mileage limit without considering exceeding the maximal maintenance time limit, because the accumulated mileage is earlier to meet the maintenance cycle than the accumulated time in general. On this basis, we modify the second part of the optimization objective function of the model. Firstly, we define a parameter $L_r$ which represents the travel mileage of the EMU circulation generated, and the parameter $r$ is the index of EMU circulation. And then, we conduct the calculation for each EMU circulation as follow: If there is $L_r \leq (1+\lambda)L_{cycle}$, we add the $(1+\lambda)L_{cycle} - L_r$ into the fitness function of the particle $F(x)$. If there is $L_r > (1+\lambda)L_{cycle}$, we add the $\beta(L_r - (1+\lambda)L_{cycle})$ into the fitness function of the particle $F(x)$, the parameter $\beta$ represents the penalty coefficient. We define a parameter $n_r$ which represents the number of EMU circulations generated, and define a 0-1 parameter $z_r$ which represents that whether the travel mileage of the EMU circulation exceeds the maintenance mileage cycle or not, if so, let $z_r = 1$, otherwise let $z_r = 0$. Therefore, by referring to the optimization objective of the model, the fitness function of the particle $F(x)$ can be described as follow:

$$F(x) = \omega_1 \sum_{i=1}^{n} \sum_{j=1, i \neq j}^{n} t_{ij} x_{ij}(1 - y_{ij}) + \omega_2 \sum_{r=1}^{n_r} \left( (1 - z_r)\left( (1+\lambda)L_{cycle} - L_r \right) + z_r \beta \left( L_r - (1+\lambda)L_{cycle} \right) \right) \tag{20}$$

Through the solution represented by the particle, we can get the connections of trains, and the maintenance arrangement is determined according to the accumulated mileage, the accumulated time and the maintenance cycle. Therefore, we define a decision variable $y_d$ for the maintenance arrangement, which represents that whether we should arrange a maintenance after assigning a train to the $d$ th position of the train connection network, if the answer is true, let $y_d = 1$, otherwise let $y_d = 0$, and this will be analyzed for the generating strategy of train connection network later on. At last, the value of $x_{ij}$ and $y_{ij}$ can be determined according to the value of $x_d$ and $y_d$. During the process of iterative optimization computation, the personal historical best solution of the particle and the global historical best solution of the particle swarm are updated by minimizing the value of $F(x)$.

### *4.2 Generating strategy of train connection network*

On the basis of particle swarm optimization algorithm, we absorb the characteristics of EMU circulation and design the generating strategy of train connection network at first. As both of the departure station and the arrival station are depot, and it is usual to arrange maintenance for EMU after undertaking an EMU circulation. Therefore, we assign a train whose departure station is adjacent to the EMU depot to the first position of the train connection network, and set this train as the first train of the corresponding EMU circulation. So maintenanceshould be arranged between the first position and the last position of the network. And then, we assign a suitable train to the following positions in sequence with meeting the constraints of train connection and maintenance cycle, and determine the maintenance arrangement reasonably. At last, we can get the train connection network with maintenance arrangement.The operation steps are shown as follows in detail:

**Step1:** According to the basic information of trains in set $V$, we select the trains whose departure station is adjacent to the depot and define a set $V_0$ consisting of these trains. Let $d=1$, and then turn to Step2.

**Step2:** If it is the process of assigning a train to the first position of the train connection network, that is $d=1$, then we choose a train $v_j$ from set $V_0$ randomly, and assign this train to the first position. And then, we remove train $v_j$ fromset $V_0$, and update the accumulated mileage $l_d$ and the accumulated time $t_d$. As it is the process of assigning a train to the first position of the train connection network, there should be amaintenance between the first position and the last position of the train connection network, in other words, a maintenanceshould be arranged after assigning a train to the $n$th position and let $y_n=1$. Therefore, the accumulated mileage $l_d$ and the accumulated time $t_d$ can be updated as $l_d=d_j^0$ and $t_d=t_j^0$ according to formula (5) and formula (6)respectively. Let $d=d+1$, and then turn to Step 3.

**Step 3:** If it is not the process of assigning a train to the first position of the train connection network, that is $2 \leq d \leq n$. Firstly, we judge that whether the arrival station of train $v_i$ assigned to the $d-1$th position of network is adjacent to the depot or not, that is to say whether there is $s_i^a \in S_{maint}$ or not. If $s_i^a \in S_{maint}$, then turn to Step 4, else if $s_i^a \notin S_{inspect}$, then turn to Step 5.

**Step 4:** If there is $s_i^a \in S_{maint}$, we choose a train $v_j$ from set $V_0$ randomly and assign it to the $d$th position of the train connection network, and remove train $v_j$ from set $V$ and $V_0$ at the same time, and update the accumulated mileage $l_d$ and the accumulated time $t_d$. If train $v_j$ meets the maintenance constraints of $l_{d-1}+d_j^0 \leq (1+\lambda)L_{cycle}$ and $t_{d-1}+t_{ij}+t_j^0 \leq (1+\lambda)T_{cycle}$, we adopt random strategy to determine whether the maintenance is arranged between train $v_i$ and train $v_j$ or not, which helps to avoid local-optimal solution as much as possible. We generate random number which belongs to the value interval $[0,1]$ at first, if the random number is less than 0.5, we arrange the maintenance between train $v_i$ and train $v_j$, and let $y_{d-1}=1$ and update the accumulated mileage $l_d=d_j^0$ and accumulated time $t_d=t_j^0$, otherwise we let $l_d=l_{d-1}+d_j^0$ and $t_d=t_{d-1}+t_{ij}+t_j^0$. If one of the maintenance constraints is not met, the maintenance must be arranged after assigning train $v_i$ to the $d-1$th position and let $y_{d-1}=1$, the corresponding accumulated mileage $l_d$ and the accumulated time $t_d$ should be updated as $l_d=d_j^0$ and $t_d=t_j^0$. Let $d=d+1$, and then turn to Step 6.

**Step 5:** If there is $s_i^a \notin S_{inspect}$, we intend to assign a train whose arrival station is not adjacent to the EMU depot to the $d$th position of the train connection network preferentially. During this process, we define two temporary sets $V'$ and $V''$ according to the trains included in set $V$ at present. The set $V'$ includes the trains which could be connected to train $v_i$ assigned to the $d-1$ th position of the train connection network, and the arrival stations of these trains are not adjacent to the EMU depot. On the contrary, the set $V''$ includes the trains which could be connected to train $v_i$ assigned to the $d-1$ th position of the train connection network, and the arrival stations of these trains are adjacent to the EMU depot. If there is a train $v_j$ in set $V'$, and it meets the constraints $l_{d-1} + d_j^0 \leq (1+\lambda)L_{cycle}$ and $t_{d-1} + t_j^0 \leq (1+\lambda)T_{cycle}$, we assign train $v_j$ to the $d$ th position of the train connection network. If there is no one train which meets the two constraints in set $V'$, we choose a train $v_j$ from the set $V''$ randomly and assign it to the $d$ th position of the train connection network. And then, we remove train $v_j$ from set $V$, and update the accumulated mileage and the accumulated time with $l_d = l_{d-1} + d_j^0$ and $t_d = t_{d-1} + t_{ij} + t_j^0$. If the set $V'$ is empty or there is no one train which meets the two constraints in set $V'$, and the set $V''$ is empty at the same time, it indicates that the train connection of the corresponding EMU circulation is completed. Therefore, we choose a train $v_j$ from the set $V_0$ randomly and assign it to the $d$ th position of the train connection network, which indicates that it is the first train of another EMU circulation, then let $y_{d-1} = 1$, update the accumulated mileage and the accumulated time with $l_d = d_j^0$ and $t_d = t_j^0$, and remove train $v_j$ from set $V$ and set $V_0$ at the same time. Let $d = d+1$, and then turn to Step 6.

**Step 6:** If there is $d > n$, it indicates that each position of the train connection network has been assigned a train, then we stop computing and output the solution of train connection network. If there is $d \leq n$, we turn to Step 3 and continue to assign a suitable train for the next position of the train connection network with the same strategy.

*4.3 Main steps of the solution algorithm*

Based on the application idea of particle swarm optimization algorithm and the generating strategy of train connection network, we optimize the solution of EMU circulation planning, and the main steps are as follows:

**Step 1:** Initializing for basic data. Input the basic information of trains, and set the maintenanceparameters and the parameters of particle swarm optimization algorithm. Then turn to Step 2.

**Step 2:** Generating for initial solution. For each particle $m$, we generate the initial train connection network with maintenance arrangement by adopting the generating strategy of train connection network presented in section 5.2, based on which we can get the initial solution of EMU circulation plan. Then we calculate the value of the fitness function $F(x)$ for each particle, and initialize the personal historical best solution for each particle and the global historical best solution for the particle swarm. Let the parameters $k=1$, $m=1$, $d=1$, and then turn to Step 3.

**Step 3:** The $d$ th dimension of the $m$ th particle indicates the $d$ th position of the train connection network represented by the $m$ th particle. Firstly, we update the value of inertia weight $w$ according to formula (19), update the velocity of the particle $v_{md}^k$ according to formula (17), and then update the train assigned to the $d$ th position of train connection network for the $m$ th particle according to formula (17). If the train $v_j$ gotten by the updating of particle swarm optimization algorithm meets the connection constraints, we assign train $v_j$ to the $d$ th position of the train connection network. If the train $v_j$ gotten by the updating of particle swarm optimization algorithm does not meet the connection constraints, then we adopt the generating strategy of train connection network presented in section 5.2 to assign a train $v'_j$ to the $d$ th position of the train connection network. According to the maintenance cycle, we determine the maintenance arrangement. And then, we update the value of $x_d$ and $y_{d-1}$, and update the accumulated mileage $l_d$ and the accumulated time $t_d$ according to the corresponding formulas. Let $d = d+1$, and then turn to Step 4.

**Step 4:** If $d \leq n$, turn to Step 3, and continue to update the train in the next position. If $d > n$, it indicates that each position of the train connection network has been assigned a train, in other words, the solution computation of the $m$ th particle has been completed. And then, we calculate the value of the fitness function $F(x)$ for the $m$ th particle according to formula (20). Let $m = m+1$, and then turn to Step 5.

**Step 5:** If $m \leq N_m$, Let $d = 1$, turn to Step 3, and continue to compute the solution represented by the next particle. If $m > N_m$, it indicates that the solution computation of all particles has been completed. And then, we update the personal historical best solution of each particle and the global historical best solution of the particle swarm according to the value of the fitness function $F(x)$. Let $k = k+1$, and then turn to Step 6.

**Step 6:** If $k \leq k_{\max}$, let $m = 1$ and $d = 1$, and then turn to Step 3. If $k > k_{\max}$, it indicates that the whole computation for optimizing EMU circulation plan has been completed, and then turn to Step 7.

**Step 7:** Stop computing and output the near-optimal scheme of EMU circulation plan.

Based on the main solving steps, we can draw a flow chart for the steps of solving EMU circulation plan, which is shown as follow:

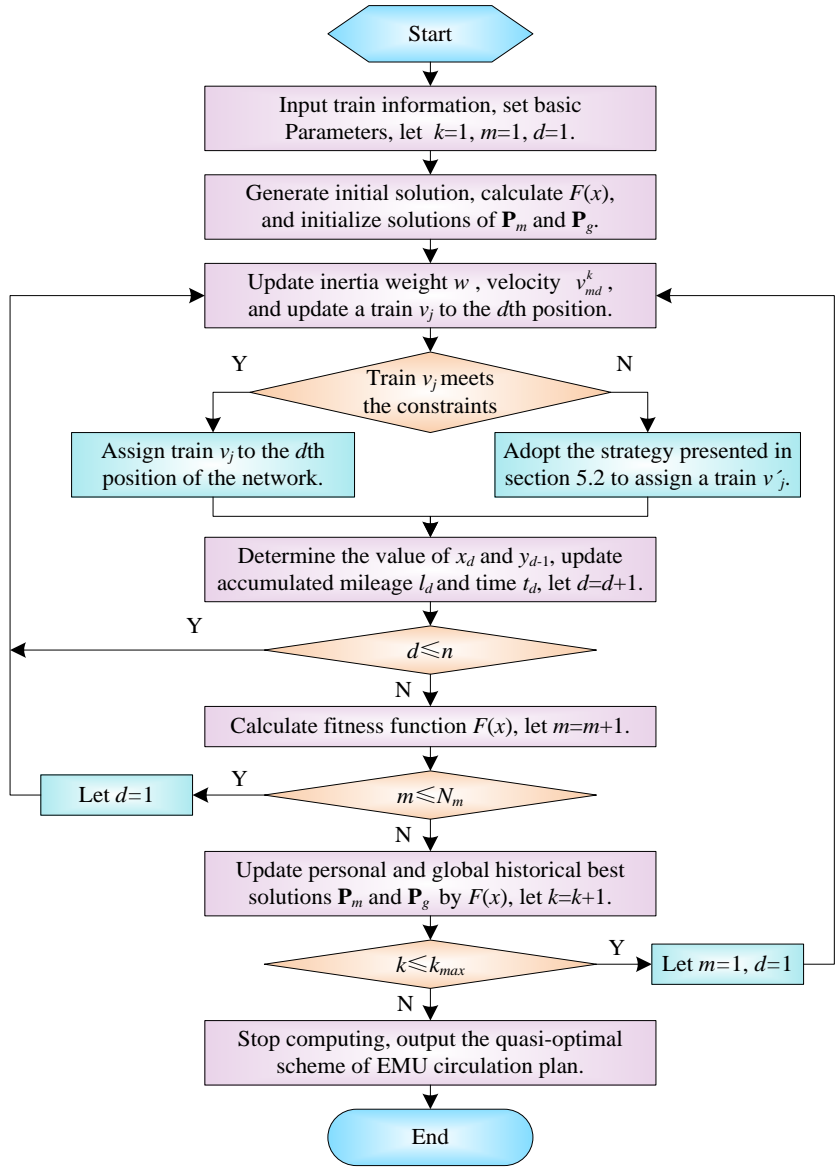

**Fig.3.** The flow chart of solving EMU circulation plan

## 5. Conclusions

In this paper, we focus on the problem of EMU circulation planning with single-type and single-depot, and we present a 0-1 integer programming mathematical model and a heuristic algorithm based on PSO for optimizing EMU circulation plan. The model is formulated on the basis of train connection network establishing, and its objectives are to minimize the amount of EMUs in service and to reduce the EMU maintenance frequency with meeting the maintenance constraints primarily.

However, this problem is a complicated one and it is affected by many other factors as well, such as the crew scheduling problem which is not considered in this paper. In future research, the problem of integrated optimization for EMU circulation planning and crew scheduling should be tackled, which helps get a more pragmatic plan. Moreover, in order to get a holistic optimization scheme for the high-speed railway system, it is necessary to integrate the EMU circulation plan, EMU assignment plan and EMU maintenance plan.

# References

Alfieri, A., Groot, R., Kroon, L., Schrijver, A., 2006. Efficient circulation of railway rolling stock. Transportation Science 40(30), 378-391.

Canca, D., Sabido, M., Barrena, E., 2014. A rolling stock circulation model for railway rapid transit systems. Transportation Research Procedia 3, 680-689.

Ceder, A. A., 2011. Optimal multi-vehicle type transit timetabling and vehicle scheduling. Procedia-Social and Behavioral Sciences 20, 19-30.

Cordeau, J. F., Soumis, F., Desrosiers, J., 2000. A Benders decomposition approach for the locomotive and car assignment problem. Transportation Science. 34, 133-149.

Cordeau, J. F., Soumis, F., Desrosiers, J., 2001. Simultaneous assignment of locomotives and cars to passenger trains. Operations Research 49(4), 531-548.

Dauzère-Pérès, S., Almeida, D. D., Guyon, O., Benhizia, F., 2015. A Lagrangian heuristic framework for a real-life integrated planning problem of railway transportation resources. Transportation Research Part B 74, 138-150.

Fioole, P. J., Kroon, L., Maróti, G., Schrijver, A., 2006. A rolling stock circulation model for combining and splitting of passenger trains. European Journal of Operational Research 174 (2), 1281-1297.

Giacco, G. L., D'Ariano, A., Pacciarelli, D., 2014. Rolling stock rostering optimization under maintenance constraints. Journal of Intelligent Transportation Systems 18(1), 95-105.

Giacco, G. L., Carillo, D., D'Ariano, A., Pacciarelli, D., Marín, Á. G., 2014. Short-term rail rolling stock rostering and maintenance scheduling. Transportation Research Procedia 3, 651 - 659.

Ho, T.K., Tsang C.W., Ip, K.H., Kwan, K.S., 2012. Train service timetabling in railway open markets by particle swarm optimization. Expert Systems with Applications 39, 861-868.

Hong, S. P., Kim, K. M., Lee, K, Park, B. H., 2009. A pragmatic algorithm for the train-set routing: The case of Korea high-speed railway. Omega 37(3), 637-645.

Kara, I., Bektas, T., 2006. Integer linear programming formulations of multiple salesman problems and its variation. European Journal of Operational Research 174(3), 1449-1458.

Lai, Y. C., Fan, D. C., Huang, K. L., 2015. Optimizing rolling stock assignment and maintenance plan for passenger railway operations. Computers & Industrial Engineering 85, 284-295.

Li, H., Han, B. M., Zhang, Q., Guo, R., 2013. Research on optimization model and algorithm of EMU circulation plan. Journal of the China Railway Society 35(3), 1-8.

Li, J., Lin, B. L., Wang, Z. K., Chen, L., Wang, J. X., 2016. A pragmatic optimization method for motor train set assignment and maintenance scheduling problem. Discrete Dynamics in Nature and Society 2016, 1-13.

Lingaya, N., Cordeau, J. F., Desaulniers, G., Desrosiers, J., Soumis, F., 2002. Operational car assignment at VIA rail Canada. Transportation Research Part B 36. 755-778.

Liu, Y. M., Niu, B., 2013. A new particle swarm algorithm theory and practice. Science Press, Beijing, China.

Miao, J. R., Wang, Y., Yang, X. Z., 2010. Research on the optimization of EMU circulation based on optimized connecting network. Journal of the China Railway Society 32(2), 1-7.

Nikolić, M., Teodorović, D., 2015. Vehicle rerouting in the case of unexpectedly high demand in distribution systems. Transportation Research Part C 55, 535-545.

Noori, S., Ghannadpour, S. F., 2012. Locomotive assignment problem with train precedence using genetic algorithm. Journal of Industrial Engineering International 8(1), 1-13.

Peeters, M., Kroon, L., 2008. Circulation of railway rolling stock: a branch-and-price approach. Computers & Operations Research 35(2), 538-556.

Rouillon, S., Desaulniers, G., Soumis, F., 2006. An extended branch-and-bound method for locomotive assignment. Transportation Research Part B 40(5), 404-423.

Salazar-González, J. J., 2014. Approaches to solve the fleet-assignment, aircraft-routing, crew-pairing and crew-rostering problems of a regional carrier. Omega 43, 71-82.

Samà, M., D'Ariano, A., D'Ariano, P., Pacciarelli, D., 2014. Optimal aircraft scheduling and routing at a terminal control area during disturbances. Transportation Research Part C 47, 1-85.

Sun, Y., Cao, C., Wu, C., 2014. Multi-objective optimization of train routing problem combined with train scheduling on a high-speed railway network. Transportation Research Part C 44, 1-20.

Thorlacius, P., Larsen, J., Laumanns, M., 2015. An integrated rolling stock planning model for the Copenhagen suburban passenger railway. Journal of Rail Transport Planning & Management 5(4), 240-262.

Wang, Z. K., 2012. Research on the optimization of operation and maintenance schemes of EMUs. Doctoral dissertation, China Academy of Railway Sciences, Beijing, China.

Wang, J. Q., 2014. LINGO-based modeling methods for the traveling salesman problem. Computer Engineering & Science 36(5), 947-950.

Zhang, C. C., Hua, W., Hua, W. J., 2010. Research on EMU scheduling under constraint of kilometrage and time for scheduled inspection and maintenance. Journal of the China Railway Society 32(3), 16-19.

Zhang, Y., Wang, S., Ji, G., 2015. A comprehensive survey on particle swarm optimization algorithm and its applications. Mathematical Problems in Engineering 2015, 1-38.

Zhao, P., Yang, H., Hu, A. Z., 1997. Research on usage of high-speed passenger trains on uncertain railroad region. Journal of the China Railway Society 19(2), 15-19.